\documentclass{itor}

\usepackage{natbib}%
\usepackage[figuresright]{rotating}

\usepackage{url}
\usepackage{algorithmic}
\usepackage{algorithm}

\usepackage{adjustbox}
\usepackage{amssymb}
\usepackage{subcaption}
\usepackage{booktabs}
\usepackage{array}

\usepackage{amsmath,amsthm}

\theoremstyle{definition}

\theoremstyle{remark}

\begin{document}

\title{A Variable Neighborhood Search for Flying Sidekick Traveling Salesman Problem}

\author[Running Author]{J\'{u}lia C\'{a}ria de Freitas\affmark{a} and  Puca Huachi Vaz Penna\affmark{a}}

\affil{\affmark{a}Departamento de Computa\c{c}\~{a}o / Universidade Federal de Ouro Preto, Campus Morro do Cruzeiro, Ouro Preto, Brazil}
\email{juliacaria@decom.ufop.br [Freitas]; puca@iceb.ufop.br [Penna]}



\begin{abstract}
An innovative model of parcel distribution is emerging from the accelerated evolution of drones and the effort of logistic companies to proceed faster deliveries at a reduced cost. This new modality originated the Flying Sidekick Traveling Salesman Problem (FSTSP) in which customers are served either by a truck or a drone. Additionally, this variant of the Traveling Salesman Problem (TSP) presents several new restrictions concerning the drone such as endurance and payload capacity.
This work proposes a hybrid heuristic that the initial solution is created from the optimal TSP solution reached by a Mixed-Integer Programming (MIP) solver.
Next, an implementation of the General Variable Neighborhood Search is used to obtain the delivery routes of truck and drone.
Computational experiments show the potential of the algorithm to improve the total delivery time up to 67.79\%. New best-known solutions (BKS) are established for all FSTSP instances that results are reported in the literature.  
Furthermore, a new set of instances based on well-known TSPLIB instances is provided.
\end{abstract}


\maketitle

\section{Introduction}\label{section:Introduction}
According to 2017 Global Online Consumer Report \citep{International2017} for 43\% of the customers, one of the most important characteristic when deciding where to buy is the enhanced delivery options. This characteristic is because the millenniums - people born between 1980 and 2000 - have a higher demand for instant satisfaction than earlier generations. Even though they are increasingly comfortable with buying products online, they are more likely to visit shops to get the product right away, rather than await delivery. 
For this reason, companies need to continually create ingenious ways to shorten delivery times and therefore satisfy the needs of customers. 
One of the most innovative announcements envisioning the improvement of delivery process occurred in December 2013 by Amazon's CEO Jeff Bezos, who declared that one of the biggest e-commerce company was testing drone parcel delivery. The project has described a drone (generically known as Unmanned Aerial Vehicles -- UAV) that departs from a warehouse loaded with a parcel, travels to the customer's location where it drops the container, and then returns to the warehouse. This operation demands no human intervention or guidance. 
Since Amazon's announcement, many companies started their drone delivery projects. Google calls Project Wing the research team responsible for developing technologies to make this drone delivery possible. Also in 2013 DHL launched a project called Parcelcopter that in 2016 accomplished more than 100 successful deliveries, including deliveries to remote villages within 30 minutes, which is faster than transporting them across steep terrains in a car \citep{BURGESS2016}. Another successful company is JD.com, China's second-largest e-commerce, that is developing a drone capable of delivering packages weighing as much as one metric ton throughout rural areas of the country flying a total of 100km before recharging \citep{Meredith2017}. Moreover, the automobile company Mercedes-Benz teamed up with the drone logistic company Matternet to start a project based on a van drone delivery concept. The project is a combination of work between van and drone in which the drone reads the destination information using a QR code on the package; then it proceeds to deliver the goods with a speed of up to 70km/h and endurance of 20km \citep{Etherington2017}.

Even though drones are fast, their payload capacity is minimal, and they are restricted by a limited endurance, this means that after each visit the drone has to return to the depot, which is not efficient. On the other hand, trucks are heavy and slow, but they can carry numerous parcels. This information is summarized in Table\ref{table:features} defined by \cite{Agatz2015}. 

\begin{table}[!htb]
\centering
\caption{Complementary features truck and drone}
\begin{tabular}{@{}lllll@{}}
\hline
      & speed          & weight         & capacity          & range \\ \hline
drone & \textbf{high}           & \textbf{light} & one               & short \\ 
truck & low   & heavy          & \textbf{many}     & \textbf{long}  \\ \hline
\end{tabular}
\label{table:features}
\end{table}

The high speed of the drone and the large capacity of the truck are complementary features of each vehicle that can make the delivery process more efficient and, possibly, at a lower cost. A possible form to take advantage of these attributes is combining the work of drone and truck to attend all customers. Hereafter is a brief description of the delivery model assumed in this work which was introduced by \cite{MurrayChu2015} known as Flying Sidekick Traveling Salesman Problem (FSTSP). Primarily, the drone is launched from the truck, then proceeds to deliver goods to a customer and finally joins back the truck in a third location. While the drone goes to a customer, the truck can travel to deliver parcels to other customers. Moreover, the truck has to travel to the returning customer before finishes the endurance of the drone.

The contributions of this work concern in combining the best features of each vehicle to build a good truck and drone delivery route to parcel delivery. 
The proposed heuristic addresses two TSP variants using drones and trucks working collaboratively, the problem presented by \cite{MurrayChu2015} and the one of \cite{Agatz2015}.
Furthermore, a new set of instances based on TSPLIB instances for Traveling Salesman Problem (TSP) is advised. 
Overall, the results obtained by the algorithm demonstrate the effectiveness of combining truck and drone for last mile parcel delivery.

\section{Related Works}\label{section:literature_review}
The use of Unmanned Aerial Vehicles (UAV) for military and surveillance purposes is broad, and exploration has been active for decades now. Nevertheless, in the last couple of years, UAV delivery research is growing, and logistics companies are investing time and money in the area. 

Here, we present some up-to-date papers dedicated to UAV delivery. The reader is refereed to \cite{Puri2005} and \cite{Mathe2015} for a comprehensive literature review on traffic surveillance and general UAV topic, respectively. For an up-to-date review related to UAV routing and trajectory, optimization see \cite{Coutinho2017}.

The Flying Sidekick Traveling Salesman Problem (FSTSP), formulated by \cite{MurrayChu2015}, describes a combination of a single truck and a single drone working in tandem. In the same paper, the authors define another issue, the Parallel Drone Scheduling TSP (PDSTSP) that is complementary to the FSTSP. 
In the PDSTSP the location of most customers is near the depot. Thus, the drone is within flight range to attend the nearest customer, and only some customers are at a higher distance, that the truck can independently serve. 
Murray and Chu proposed a heuristic for both problems, and for the FSTSP they also introduced a Mixed-Integer Linear Programming (MILP) formulation. 

\cite{Ponza2016} based his dissertation on the FSTSP proposed by Murray and Chu. He proposed a slightly different mathematical formulation to the problem and presented an analysis of several heuristics that could be used to resolve the problem. Moreover, he introduced a new set of instances to literature.

The TSP-D formulation proposed by \cite{Agatz2015} assumes that both truck and drone travel on the same road network. Furthermore, the drone may launch and return to the same location. 
The work was extended in \cite{Bouman2017} presenting two different approaches to solve the problem. The dynamic programming approach consists of finding a solution to the original TSP, subsequently combining this truck path with drone nodes to obtain a route where the truck nodes are covered using the shortest path algorithm. Furthermore, for the purpose to overcome the run time of dynamic programming, they also proposed an $A^*$ algorithm which was faster for larger instances in all cases. 

\cite{Ha2015} describe a variant of the FSTSP. They proposed a MILP formulation and two methods to solve the problem, a GRASP based heuristic, and a TSP-LS. Computational experiments performed in instances with up to 100 customers show that GRASP outperforms TSP-LS regarding solution quality, but regarding running time, GRASP presents a greater computational time.

The work of \cite{Ferrandez2016} envisions reduce delivery time and energy of the work in collaboration of truck and drone (in-tandem system) with the traditional truck delivery (stand-alone). The authors proposed a genetic algorithm (GA) to compute the stand-alone delivery, i.e., the TSP truck route. Furthermore, a K-means algorithm is used to find a solution that combines the truck route with drone nodes.

\cite{Mathew2015} discuss two different problems in their work. First, a Heterogeneous Delivery Problem (HDP) as a discrete optimal path planning problem that seeks to minimize the total cost of deliveries in routes computed for a truck and a drone. 
To solve the HDP, they devise a reduction to obtain a Generalized Traveling Salesman Problem (GTSP) and be able to use already available solvers. This paper uses three different edges cost: one for drone flight edges, one for truck street edges and one for when the truck transports the drone. The second problem presented is called Multiple Warehouse Delivery Problem (MWDP). In MWDP only drones perform the delivery of parcels,  again they use reduction to GTSP to solve this problem with known solvers.

\cite{Ulmer2017} present a different approach to the problem, a Same Day Delivery Problem (SSDP). For each ordering customer, the provider must decide if the order will be performed in the same day or not and for which vehicle, the truck or the drone. 
The authors define a Markov decision process model to describe the problem and use a policy function approximation (PFA) to determine the best values for a parameterized policy. The parameter is a threshold of travel distance used to define which vehicle will attend each customer. Trucks preferably serve the customers in the zone of the threshold, and drones preferably attend customers with further distances. 
The results show that the districting by the proposed PFA is beneficial once it is possible to increase the number of services by the fleet significantly.

\cite{Wang2017} study the Vehicle Routing Problem with Drones (VRPD) from a worst-case point of view. The paper describes several theorems formulation for the vehicle routing problem with drones and represents bounds on maximal savings to the companies. \cite{Poikonen2017} enlarger the description of the theorems comparing different drone configurations in the delivery process to determine the maximum benefit. For example, the trade-off between speed and the number of drones, i.e., they compare what is better, a more substantial number of slower drones or a smaller number of faster drones. 

\cite{DiPugliaPugliese2017} provide an Integer Programming formulation for the Vehicle Drone Routing Problem with Time Windows (VDRPTW) with the objective function of minimizing the total transportation cost. In this work, the drone has a waiting time limit of the truck after performing a delivery, which prevents it from being idle. Additionally, after delivery, the drone must return to the truck it was launched that is located at another customer.

The problem introduced by \cite{Daknama2017} uses several trucks and drones traveling along a route. Trucks follow Manhattan metrics while drones follow Euclidean metrics. A local search was implemented in a framework called JAMES, customized with appropriate neighborhoods.

The primary goal of \cite{Goodchild2017}'s paper is to determine whether or not drone technology in the logistics industry would have a positive environmental impact regarding vehicle-miles traveled (VMT) and carbon dioxide ($CO_2$) emissions. Thus, the work presents models with different scenarios in which trucks and drones originating from a central depot deliver parcels to recipient addresses in circular service zones. 
The results suggest that a system would perform best with drones serving nearby addresses and trucks delivering to the farther ones.

Another variant for the problem is the Drone Delivery Problem (DDP) by \cite{Dorling2017} where only drones perform the delivery. 
The distribution center is the return point of routes, where the delivery man changes the battery of the UAV and loads it with another parcel. Moreover, drones can visit a location only once, and perform multiple trips per route. The paper provides a MILP implementation for small instances and for large ones a SA heuristic was proposed to find suboptimal solutions to the DDP within a limited run time.

\cite{Vorotnikov2017} also consider the drone as the only vehicle. The addressed problem determines the drone route by solving the TSP with three approaches: the Monte Carlo method, the method of reduction of rows and columns and the averaged coefficient method. Research results show that reduction of rows and columns presents the best cost. However, according to the time criterion for the solution, the latter has a significant advantage, which is especially noticeable with an increase in the number of objects.

\cite{Othman2017} introduce four different variants for DDPs. In the first, the drone launches from the truck, and while it is delivering a parcel, the truck is allowed to wait for the drone to return at the previous rendezvous point whereas in the second description the truck is not allowed to wait at the previous rendezvous point. 
The others definitions are similar, in the third, the truck is allowed to wait for the drone to return at the previous rendezvous point and the drone can be transported by the truck before proceeding to its next customer. 
Finally, in the last one, the truck is not allowed to wait for the returning of the drone at the previous rendezvous point.
They modeled this problem as a problem of finding a special type of a path in a graph of a special structure, and proposed a polynomial-time approximation algorithm for the graph problem in metric graphs.

 \cite{Coelho2017} present an approach using the concept of smart cities to introduce a multi-objective green UAV routing problem. The problem considers a dynamic scenario in which new orders may arrive at any moment. The scenario is composed of airspace divided into two layers: a lower layer, located at low altitudes with lower speed drones and an upper layer, located at high altitudes composed of faster drones. A Multi-Objective Smart Pool Search (MOSPOOLS) Matheuristic is used to obtain a solution to the proposed problem.

The use of drones in a humanitarian manner can have a significant role, especially in post-disaster operations. Whereas it is possible to take advantage of the favorable characteristics of drones, such as the high speed to deliver urgently needed small items in locations with difficult access.

The work of \cite{Scott2017} is concerned with creating drone delivery models for health-care to locate a warehouse with supplies and drone nests to complete final delivery. They developed two models to the problem; the first one has an objective to minimize the total weighted delivery time, i.e., road covered by the truck plus air covered by the drone. The second model intends to minimize the maximum weighted time for truck/drone delivery subject to a budget and drone travel distance constraint.

Table \ref{table:papers} presents an overview of the works related to drone delivery. The authors and year of the publication are in the first column. The second column presents the models proposed. The columns three and four indicate how many vehicles each problem uses. The fifth column states if the node the drone launch from can be the same it returns to. 
The sixth column indicates if the vehicles follow the same network. Finally, the last column shows the method each author used to address the delivery problem.
The column "Problem" presents the acronym adopted in the respective paper or the type they belong to if no definition is provided. Different acronyms have the same description, as can be seen in \cite{Ha2015} and \cite{Ferrandez2016}. 

The table summarizes the differences and similarity between the formulations using delivery with drones. A distinct number of drone and truck is used for each work as well as the road network the vehicles travels in. This work deals two variants, the one introduced by \cite{MurrayChu2015} and the one proposed by \cite{Agatz2015}.

\begin{table}[!htb]
\caption{Related Works and characteristics. }
\label{table:papers}
\begin{adjustbox}{max width=\textwidth}
\begin{tabular}{c|c|c|c|c|c|c}
\hline

      &         & \#         & \#          &  Launch             & Drone network          & Proposed \\
\multicolumn{1}{c|}{Authors} & Problem & Truck     &  Drone    &  $\neq$               & $\neq$                & Solution \\
      &           &            &              &  Return            & Truck network         & \\ \hline
      
\cite{MurrayChu2015} & FSTSP &1 & 1 & \checkmark & \checkmark   & MILP formulation and\\
                     & PDSTSP &1 & 1 & \checkmark & \checkmark   & Heuristic approach\\ \hline
\cite{Ponza2016}      & FSTSP &1 & 1 & \checkmark &   &  SA \\ \hline
&&&&&& Kruskal's MST \\  
\cite{Agatz2015}     & TSP-D &1 & 1 &               &           & Heuristic and DP\\
                                                     &&&&&& MILP formulation \\ \hline
\cite{Bouman2017} & TSP-D &1 & 1 &               &           & Dynamic Programming \\
                                                      &&&&&& A* \\ \hline
                 &&&&&&  GRASP \\ 
\cite{Ha2015}   & TSP-D &1 & 1    & \checkmark  & \checkmark  & TSP-LS \\
                &&&&&& MILP formulation \\\hline
\cite{Ferrandez2016} &TSP-D&1 & m & \checkmark & \checkmark &   GA and K-means\\ \hline

\cite{Mathew2015} & HDP & 1 & 1 & & \checkmark &  Reduction to GTSP \\
                  & MTSP& 0 & m & & \checkmark &  Reduction to TSP \\ \hline
\cite{Ulmer2017} & SDDPHF & n & m & & &  Parametric policy\\
                                   &&&&&&  function approximation \\ \hline

\cite{Wang2017} & VRPD & n & m & & &  Theorems with \\ 
                                     &&&&&& worst case scenarios \\\hline
\cite{Poikonen2017} & VRPD & n & m & & &  Theorems with \\ 
                                     &&&&&& worst case scenarios \\\hline
\cite{Daknama2017} & VRD & n & m & \checkmark & \checkmark &  Local Search\\ \hline
\cite{Dorling2017} & DDP & 0 & 1 & & &  SA \\ 
                                   &&&&&& MILP formulation \\ \hline
                                         &&&&&& Monte Carlo  \\                              
 \cite{Vorotnikov2017}                   & TSP & 0 & 1 & & &                Reduction of rows  \\ 
                                  &&&&&& and columns \\ 
                                                                     &&&&&& MILP formulation \\ \hline
 & ALSDP & 1 & 1 & & &  2-approximation \\
\cite{Othman2017}                  & NW-ALSDP & 1 & 1 &\checkmark & &  Min cost bipartite method \\
                  & LSDP & 1 & 1 & & &  2-approximation \\
                  & NW-LSDP     & 1 & 1 &\checkmark & &   Min cost bipartite method \\ \hline
    
\cite{DiPugliaPugliese2017} & VDRPTW & n & m & \checkmark & \checkmark &  IP formulation \\  \hline
This work & FSTSP & 1 & 1 & \checkmark & \checkmark &  GVNS \\   
		  & TSP-D & 1 & 1 &  &  &   \\ 
\hline

\end{tabular}
\end{adjustbox}
\end{table}

\section{Problem Description}\label{section:description}

This work assigns a variant of the classical Traveling Salesman Problem (TSP) called Flying Sidekick Traveling Salesman Problem (FSTSP) first introduced by \cite{MurrayChu2015}.

The FSTSP, also considered an NP-Hard problem, describes a drone and a truck that collaboratively delivery parcels to a set of customers. One drone and one truck execute the deliveries. The drone can launch and return to the truck at any customer. However, the drone can only visit eligible customers, i.e., customers that the parcel does not exceed the vehicle payload capacity and respects the endurance of the vehicle to complete the trip. Additionally, while the drone goes to a node, the truck can independently travel to deliver parcels to others customers, considering that the vehicle proceeds to the return node before finishes the endurance of the drone.

For better understanding, we define some parameter notation to represent the sets used in the problem.
Let $C = \{1, ..., c\}$ denotes the set of all customers and let the subset $C^{\prime} \subseteq C$ represents the eligible customers to the drone. The distribution center (depot) is denoted by node 0. Thus $N = \{0, 1, ..., c\}$ represents all nodes in the network. 

The drone is allowed to fly in a straight line ignoring the road restrictions; however, the truck must follow the road network. Additionally, vehicles may have different speed. Thus, we consider the travel time instead of the travel distance to respect the traffic characteristics. The truck travel time from node $i$ to node $j$ is given by $\tau_{ij}$. Analogously, $\tau_{ij}^\prime$ represents the travel time of the drone.
Moreover, we consider service time, i.e., the time required before the drone launches to change a battery and load the vehicle with a parcel, represented by $s^l$. $s^r$ is the time necessary to recover the drone after a delivery. 
The drone has a limited flight endurance defined by parameter~$e$.

The objective is to find the minimum cost route that serves all customer locations by either the truck or drone. Both the departure and return to the depot can occur in tandem or independently. However, while traveling in tandem, the truck must transport the drone to save battery. As the drone has unitary capacity, it has to pick up a new parcel at the truck after each delivery. 

\begin{figure}[!hpt]
	\centering
    \includegraphics[width=\textwidth]{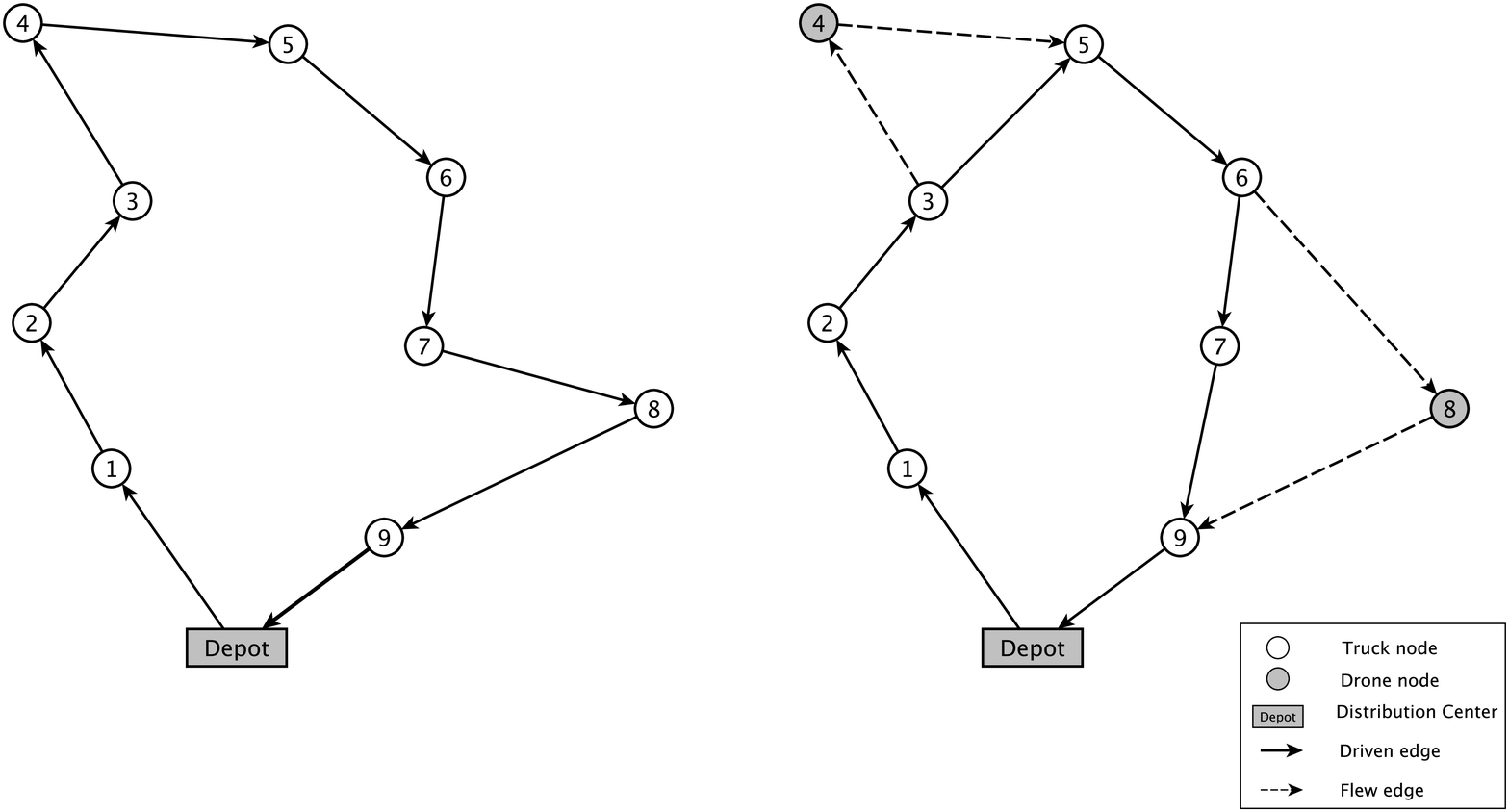}
    \caption{Flying Sidekick Traveling Salesman Problem.}
    \label{fig:routes}
\end{figure}

Figure \ref{fig:routes} represents a delivery route where nine customer locations are needing to be attended. It is reasonable to observe that serving two customers (customers 4 and 8) with the drone instead of the truck reduces the total distance traveled by truck. Therefore, by shortening the truck route, it is possible to reduce the overall delivery time required to serve all customers. 

As demonstrated in Figure \ref{fig:routes} a customer can be truck-only, drone-only or mixed. \textit{Mixed} are the customers where the drone launch or return (customer 3, 5, 6, 9). \textit{Drone-only} customers are the ones visited only by drones (customer 4, 8). Finally, \textit{truck-only} customers are the remaining truck customers (customer 1, 2, 7).

The drone has a trip composed of three distinct nodes: launch, visiting, and return. The launch node is where the driver can change the battery and load the drone with a parcel. A visiting node consists of a customer that is visited only by the drone. The return node describes the location where the driver recovers the drone. 
The return node can either be a customer serviced by the truck or the depot, though it can not be a customer already visited. Furthermore, in the case of the depot, the drone cannot be relaunched. All drone trips must respect the endurance flight limit, i.e., the drone has to have enough battery to visit the three nodes and, when necessary, wait for the truck in the return node. 

\begin{figure}[!ht]
  \begin{subfigure}[b]{0.48\textwidth}
    \centering
    \includegraphics[width=0.8\textwidth]{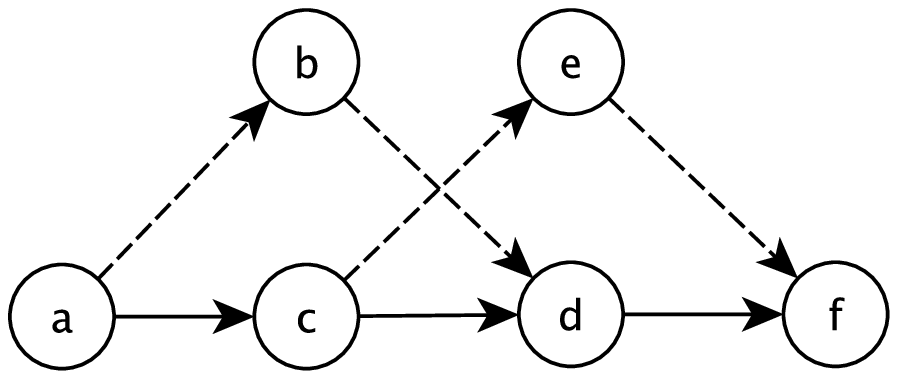}
    \caption{Prohibition 1.}
    \label{fig:prohibition1}
  \end{subfigure}
  \hfill
  \begin{subfigure}[b]{0.48\textwidth}
    \centering
    \includegraphics[width=0.8\textwidth]{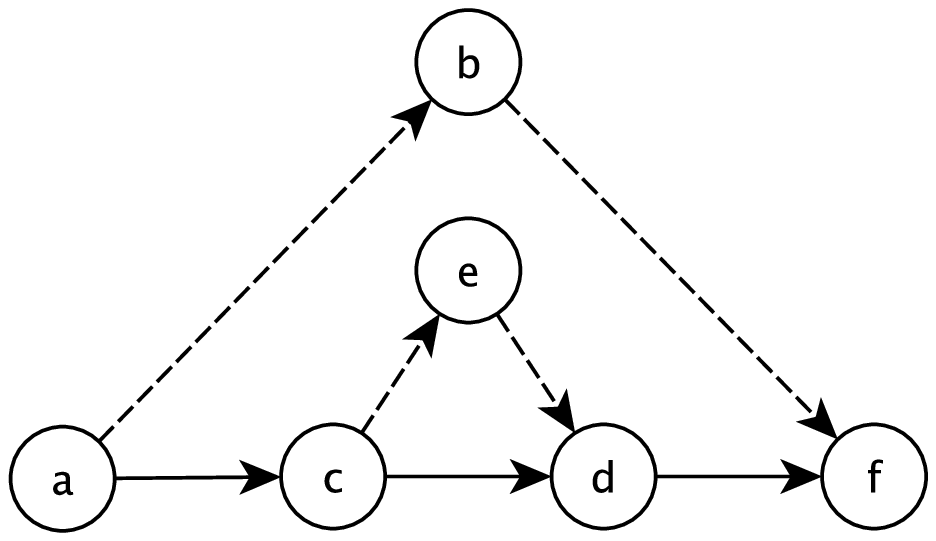}
    \caption{Prohibition 2.}
    \label{fig:prohibition2}
  \end{subfigure}
  \caption{Drone Trip Prohibitions.}
  \label{fig:prohibitions}
\end{figure}

As the FSTSP considers only one drone, two prohibited situations exist in the route, as reported by \cite{Ponza2016}. Figure~\ref{fig:prohibitions} describes those situations within continuous line representing the truck tour and dashed lines the drone tour. The first case represented in Figure \ref{fig:prohibition1} happens when a drone trip starts before the last trip finishes, i.e., a drone is launched from customer $c$ when the trip ($a, b, d$) was not finished yet. 
The other situation occurs when a drone trip starts and ends before the last trip finishes, i.e., a launch (node~$c$) and a return (node~$d$) is inside another trip of launch-return (nodes $a$ and $f$) described in Figure~\ref{fig:prohibition2}.

\section{HGVNS Algorithm}
 \label{section:method}

The algorithm framework used in this work presents an exact model to obtain an initial solution and the General Variable Neighborhood Search (GVNS) metaheuristic as an improvement heuristic.
Algorithm~\ref{algorithm:Framework} illustrates the hybrid GVNS heuristic framework, named HGVNS.
The algorithm is a simple hybrid heuristic integrating an exact method and heuristics. The HGVNS is composed of three steps, described below.
The input data is the set of customers and the cost matrices. Two different cost matrices are used: the truck and the drone travel time, $\tau$ and $\tau^{\prime}$, respectively (see Section~\ref{section:description}). 

The first step of the algorithm creates an initial solution where the truck visits all customers. Since this is a regular TSP, a well-studied problem, we use a Mixed-Integer Programming (MIP) solver to determine the optimal TSP route. 
The solver requires the truck travel time matrix $\tau$ and the entire set of customers to generate the solution (line~\ref{algorithm:F-concorde}). The solution $s^*_{TSP}$ is the optimal solution for the TSP, as it is also a feasible solution for the FSTSP, but not including any drone trip. $s^*_{TSP}$ is given to $CreateInitialSolution()$ procedure to add the drone trips to the route (line~\ref{algorithm:F-InitSol}). Finally, the initial solution $s$ is given to the GVNS procedure (line~\ref{algorithm:F-VNS}), to obtain the improved solution~$s^*$, which is then returned by HGVNS.

\begin{algorithm}[!htb]
\caption{HGVNS Framework}
\label{algorithm:Framework}
\begin{algorithmic}[1]

   \REQUIRE cost matrices ($\tau$, $\tau^{\prime}$); customers ($C$, $C^\prime$)
   \STATE $s^*_{TSP} \leftarrow $ MIP\_Solver($\tau$, $C$)\label{algorithm:F-concorde}
   \STATE $s \leftarrow$ CreateInitialSolution($s^*_{TSP}$, $\tau$, $\tau^{\prime}$, $C$, $C^\prime$) \label{algorithm:F-InitSol}
   \STATE $s^* \leftarrow$ GVNS($s$, $\tau$, $\tau^{\prime}$, $C$, $C^\prime$) \label{algorithm:F-VNS}
   \RETURN $s^*$
\end{algorithmic}
\end{algorithm}

The $CreateInitialSolution()$ procedure, based on the heuristic presented by \cite{MurrayChu2015}, removes some of the truck's customers and makes them drone customers. Algorithm \ref{algorithm:initialSolution} illustrates the procedure. The algorithm initializes variable $maxSavings$ and $truckSubRoutes$. 
Following, for each eligible drone customer $j$ is computed the savings of removing this customer from the truck route (line \ref{initialSolution:savings}). 

\begin{algorithm}[!htb]
    \caption{CreateInitialSolution}
    \label{algorithm:initialSolution}
    \begin{algorithmic}[1]
    
       \REQUIRE cost matrices ($\tau$, $\tau^{\prime}$); customers ($C$, $C^\prime$) 
       \ENSURE s
      
      \STATE $maxSavings \gets$ 0
      \STATE $truckSubRoutes \leftarrow$ \{$s^*_{TSP}\}$
      \REPEAT
        
        \FORALL{$j \in$  $C^{\prime}$ } \label{initialSolution:loop}
           \STATE \label{FSTSP:savings} $savings \leftarrow$ computeSavings($truckSubRoutes$, $j$)\label{initialSolution:savings}
            \FORALL{$subroute \in$ $truckSubRoutes$}
                \IF {(isPairedUAV($subroute$)}
                    \STATE $maxSavings \gets$ computeCostTruck($subRoute$, $j$) \label{initialSolution:truck}
                \ELSE 
                    \STATE $maxSavings \gets$ computeCostUAV($subRoute$, $j$) \label{initialSolution:drone}
                \ENDIF
            \ENDFOR
        \ENDFOR \label{initialSolution:endloop}
        
        \IF{$maxSavings > 0$}
            \STATE Update the route\label{initialSolution:update}
            \STATE Reinitialize $maxSavings \gets 0$
        \ENDIF
      \UNTIL {Stop}
      
    \end{algorithmic}
\end{algorithm}

Next, for each sub-route, it is verified if there is a drone trip (lines~\ref{initialSolution:loop}~--~\ref{initialSolution:endloop}). If the answer is positive, customer $j$ may be relocated in the truck route. If the answer is negative, customer $j$ may be visited by the drone, creating a new sub-route. Ultimately, if one of the operations above reduces the delivery time, the routes of the vehicles are updated (line \ref{initialSolution:update}).

The GVNS presented in the third step of HGVNS is defined in Algorithm~\ref{algorithm:VNS}. 
Given an initial solution (line~\ref{VNS:input}), the main operation iterates over the parameter $k$ until it reaches the stop condition, which is the cardinality of the list of neighborhoods $|\mathcal{N}|$. 
Each interaction generates a solution $s^\prime$ from a neighborhood $\mathcal{N}^{(k)} (s)$ of current solution $s$ (line~\ref{VNS:randomN}).

Next, a local search, named RVND (Algorithm~\ref{algorithm:RVND}), is performed in solution $s^\prime$ to obtain $s^{\prime\prime}$ as the local optimum (line \ref{VNS:LO}). If the cost value of the obtained solution is better than the incumbent one, $s^{\prime\prime}$ is assigned as the current solution and the search continues with $\mathcal{N}^{(k)} (s)$. Otherwise, the next neighborhood in the list $\mathcal{N}$ is executed (lines~\ref{VNS:ifs} -- \ref{VNS:ife}). Both procedures, GVNS and RVND use the same neighborhoods list~$\mathcal{N}$.

\begin{algorithm}[!htb]
    \caption{General Variable Neighborhood Search}
    \label{algorithm:VNS}
    \begin{algorithmic}[1]
    
       \REQUIRE $s$, $\tau$, $\tau^{\prime}$, $C$, $C^\prime$, $k_{max}$ \label{VNS:input}
       \ENSURE s
      \STATE Initialize Neighborhood List ($\mathcal{N}$)
      \STATE $k \gets$ 1
      
      \WHILE{$k \le k_{max}$} 
          \STATE  Generate a point $s^{\prime}$  at random from $k^{th}$ neighborhood of s ($s^{\prime} \in \mathcal{N}^{(k)}(s)$). \label{VNS:randomN}   
        \STATE  $s^{\prime\prime} \gets$ RVND($s^{\prime}$, $\tau$, $\tau^{\prime}$, $C$, $C^\prime$) \label{VNS:LO}
        
        \IF{$s^{\prime\prime} < s$}\label{VNS:ifs}
            \STATE $s \gets s^{\prime\prime}$
            \STATE $k \gets 1$
        \ELSE
            \STATE $k \gets k+1$    
        \ENDIF \label{VNS:ife}
      \ENDWHILE
    \end{algorithmic}
\end{algorithm}

The local search used by GVNS is the Randomized Variable Neighborhood Descent (RVND). The RVND, presented in Algorithm \ref{algorithm:RVND}, differs from VND by randomly choosing the next neighborhood. First, the heuristic initializes a list of neighborhoods and then shuffles it (line~\ref{RVND:neighborhood_list}). The algorithm ends when counter $k$ reaches the stop condition, which in this case is to perform all neighborhoods without improvement (line~\ref{RVND:while}). We apply the Best Improvement (BI) approach that exhaustively explores the neighborhood and returns one of the solutions with the lowest solution value (line \ref{RVND:best}).
Then, the neighborhood solution  $s^\prime$ is compared with the current solution $s$. Thus, if $s^\prime$ is better than $s$, $s^\prime$ is the new solution. The list of neighborhoods is reinitialized and shuffled again (lines~\ref{RVND:solution} -- \ref{algorithm:RVND-shuffle2}) and the counter is restarted $k$. The neighborhoods used in this method are detailed in the Sections~\ref{subsection:Reinsertion} to \ref{subsection:relocateCustomer}.

\begin{algorithm}[!htb]
    \caption{Randomized Variable Neighborhood Descent}
	\label{algorithm:RVND}
    \begin{algorithmic}[1]
    
       \REQUIRE $s$, $\tau$, $\tau^{\prime}$ $C$, $C^\prime$
       \ENSURE $s$
       
      \STATE Initialize Neighborhood List ($\mathcal{N}$) \label{RVND:neighborhood_list}
      \STATE $k \gets$ 1
      
      \WHILE{($k \le \mathcal{N}$)} \label{RVND:while}
          \STATE  Find the best neighbor $s^{\prime} \in \mathcal{N}^{(k)}(s)$ \; \label{RVND:best}
        
        \IF{$s^{\prime} < s$}\label{RVND:solution}
            \STATE $s \gets s^{\prime}$
            \STATE Reinitialize and shuffle $\mathcal{N}$ \label{algorithm:RVND-shuffle2}\;
            \STATE $k \gets 1$
        \ELSE
            \STATE $k \gets k+1$    
        \ENDIF 
      \ENDWHILE
    \end{algorithmic}
\end{algorithm}

\subsection{Neighborhoods}
\label{section:neighborhoods}

The data structure used to store the solution of both truck and drone is an array. However, the drone's array is composed of a tuple of three values, indicating the launch, visit and return node of a trip.

When a change is performed on a truck-only customer the cost of the new solution can be calculated in $\mathcal{O}(1)$, using the following formula:

\begin{equation}
\label{equation:cost}
\begin{split}
\centering
    cost = f(s) + \left(\sum_{i \in \theta^+}{\tau_{i}} - \sum_{j \in \theta^-}{\tau_{j}} \right)
\end{split}
\end{equation}

In Equation~\eqref{equation:cost}, set $\theta^-$ represents the removed edges and $\theta^+$ the reconnected edges set. Sets $\theta^-$ and $\theta^+$ have a fixed size. Considering $n$ the number of customer, $|\theta^-| << n$ and $|\theta^+| << n$ .
An example of the use of Equation~\eqref{equation:cost} is to calculate the new cost of the solution presented in Figure \ref{fig:reinsertion2}(b). The edges \{(0,5), (3,1), (6,2)\}  $\in$ $\theta^-$ were removed and edges \{(0,1), (6,5), (3,2)\} $\in$ $\theta^+$ were reconnected. Therefore, the new solution cost is represented in Equation~\ref{equation:newcost}, considering that $\tau$ is the travel time between two customers.

\begin{equation}
\label{equation:newcost}
\begin{split}
\centering
     cost = f(s) + \left( (\tau_{01} + \tau_{65} + \tau_{32}) - (\tau_{05} + \tau_{31} + \tau_{62})  \right)
\end{split}
\end{equation}

Now, if the movement is performed in mixed customers, different situations arise. However, the drone trip only modifies the solution cost when its last trip returns to the depot. Thus, even when the drone trip is modified, the new solution cost may be defined entirely by the changes in the route of the truck. 

Hereafter is presented the seven neighborhoods structures. A neighborhood movement is accepted as long as it results in a feasible solution: the move must not violate the endurance constraint, and it must not create the prohibited situations of Figure~\ref{fig:prohibitions}.

\subsection{Reinsertion}
\label{subsection:Reinsertion}
\begin{figure}[H]
    \centering
    \includegraphics[width=0.9\textwidth]{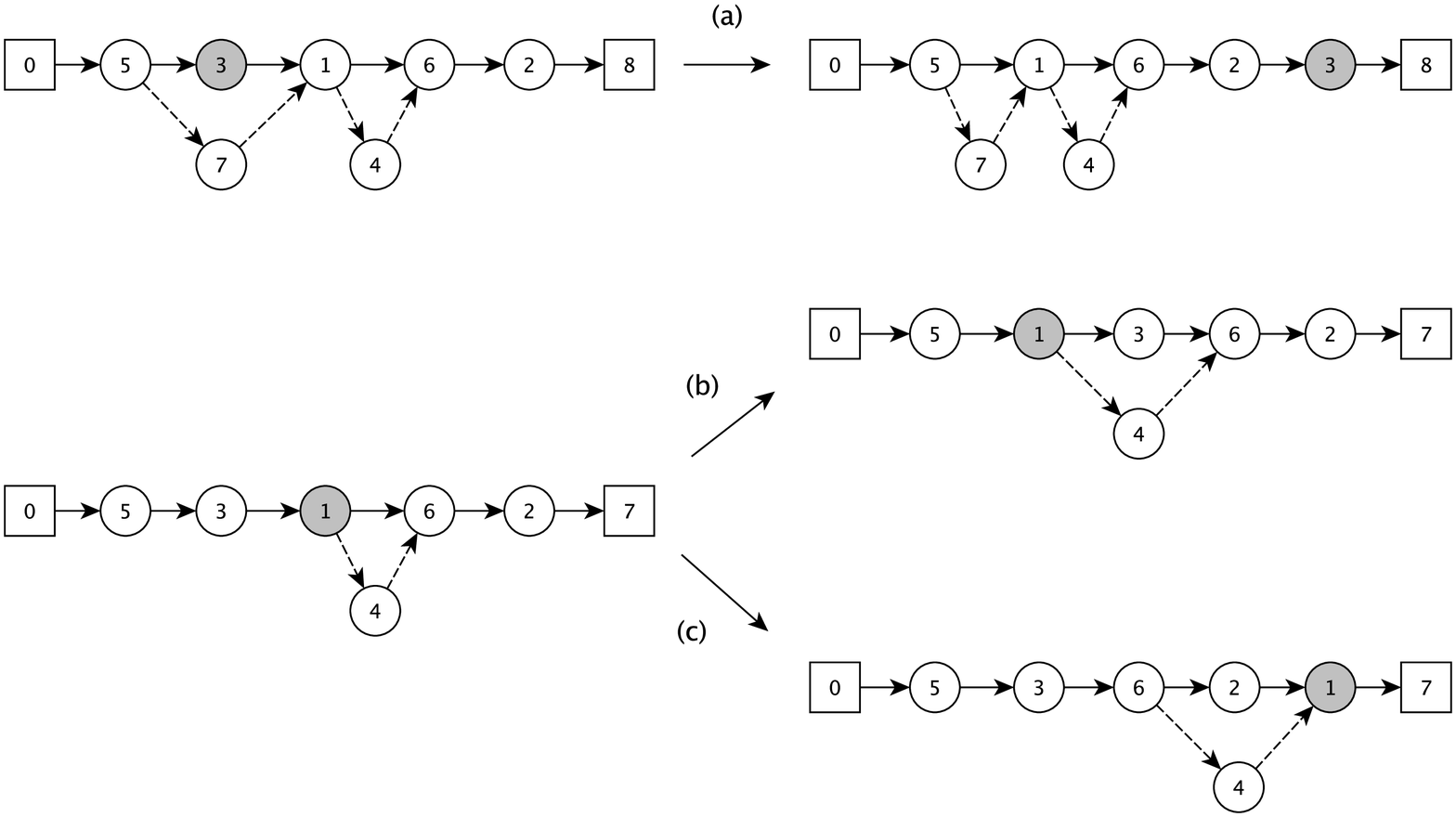}
    \caption{Move Reinsertion.}
    \label{fig:reinsertion}
\end{figure}

This neighborhood removes a customer and reinserts it in other position in a tentative solution.
Figure~\ref{fig:reinsertion}(a) illustrates a relatively straightforward move in which a truck-only customer is relocated.
The path represented by Figure \ref{fig:reinsertion}~(b) -- (c) relocates the customer 1, additionally, in both relocations the truck route during a drone trip increases, thus, the drone trip is changed. In case (b) the drone still launches from customer 1, however, in case (c) the trip is inverted, turning customer 1 the return node of a drone's trip.

\subsection{Or-opt2}
\label{subsection:or-opt2}
\begin{figure}[H]
    \centering
    \includegraphics[width=0.9\textwidth]{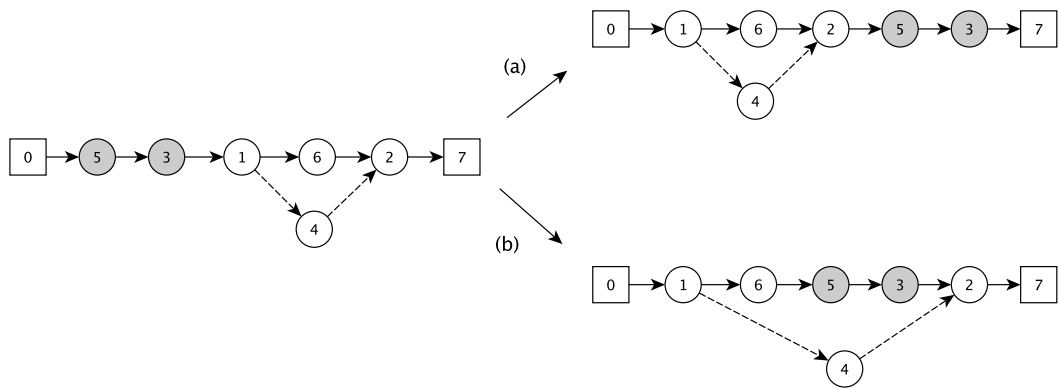}
    \caption{Move Or-opt2.}
    \label{fig:reinsertion2}
\end{figure}

    This neighborhood relocates two adjacent nodes of the truck path in an arbitrary position in a tentative solution. In Figure~\ref{fig:reinsertion2} the consecutive nodes 5 and 3 are relocated first in example (a) not affecting the route of the drone and in example (b) increasing the sub-route with a drone trip.
    
\subsection{Exchange} \label{neighborhood:exchange}
\label{subsection:exchange}
 This neighborhood is set to swap a customer with another one in a tentative solution. A violation of the prohibition \ref{fig:prohibition2} would occur if the launch and return nodes remain the same. Thus, customer 2 was chosen to be the new returning node of the drone's trip, this way, the trip \{5, 7, 1\} became \{5, 7, 2\}. 
Moreover, the second drone trip \{1, 4, 6\} continues the same, as the distances are symmetric.

\begin{figure}[H]
    \centering
    \includegraphics[width=0.9\textwidth]{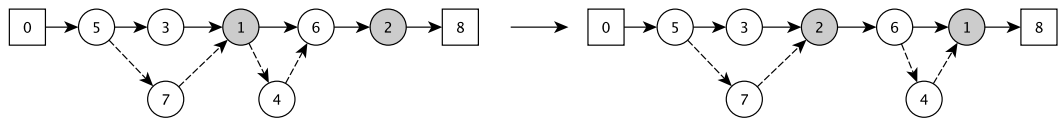}
    \caption{Move Exchange.}
    \label{fig:exchange}
\end{figure}

\subsection{Exchange(2,1)}
\label{subsection:exchange21}
  This neighborhood is set to swap two adjacent customers with another customer in a tentative solution. Following the explained in neighborhood \ref{neighborhood:exchange}, the first trip \{5, 7, 1\} must have another return node to not violate prohibition shown in Figure~\ref{fig:prohibition2}. In this, case node 2 was the one selected. Furthermore, the trip \mbox{\{6, 3, 1\}} was inverted, and the travel distance of the truck enhanced, increasing the waiting time of the drone.
 
\begin{figure}[H]
    \centering
    \includegraphics[width=0.9\textwidth]{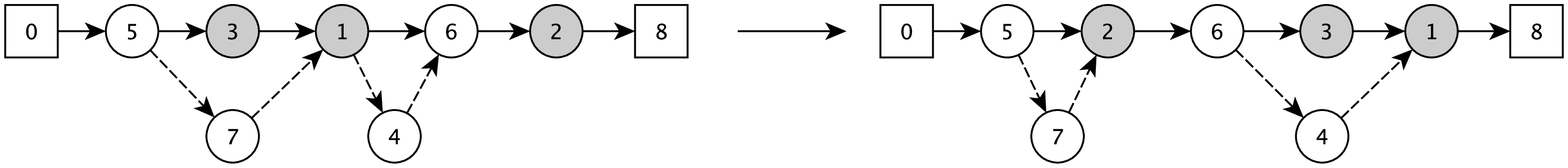}
    \caption{Move Exchange(2,1).}
    \label{fig:exchange21}
\end{figure}

\subsection{Exchange(2,2)}
\label{subsection:exchange22}

\begin{figure}[H]
    \centering
    \includegraphics[width=0.9\textwidth]{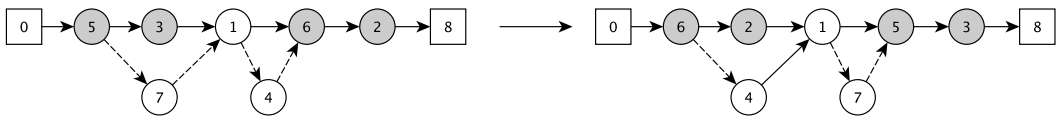}
    \caption{Move Exchange(2,2).}
    \label{fig:exchange22}
\end{figure}

This neighborhood is set to swap two adjacent customers with another two adjacent customers in a tentative solution. In the example illustrated by Figure~\ref{fig:exchange22} both trips were reversed. While the drone waiting time travel of trip \{5, 7, 1\} decreased, the waiting time of the trip \{1, 4, 6\} increased.

\subsection{2-opt}
\label{subsection:2opt}

\begin{figure}[H]
    \centering
    \includegraphics[width=0.9\textwidth]{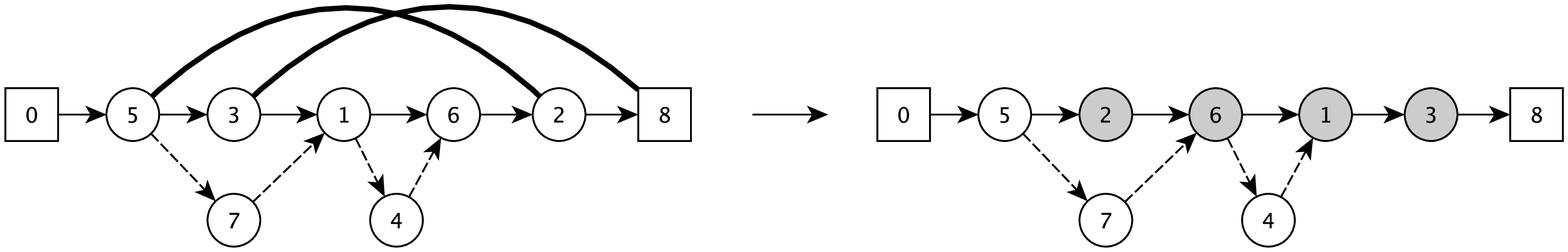}
    \caption{Move 2-opt.}
    \label{fig:2opt}
\end{figure}

This move performs a 2-opt on the solution: two edges are removed from a tentative solution and the two paths created are reconnected in the only possible way to keep a valid tour. The removed edges must be truck-only or mixed nodes, and if one of these edges are under a drone trip another return node must be selected (represented by trip \{5, 7, 1\} that changed to \{5, 7, 6\} in Figure~\ref{fig:2opt}).  Figure~\ref{fig:2opt} represents the 2-opt move where edges (5, 3) and (2, 8) are removed. Thus edges (5, 2) and (3, 8) are created, reversing the path between nodes 5 and 8. 

\subsection{Relocate Customer}
\label{subsection:relocateCustomer}

\begin{figure}[H]
    \centering
    \includegraphics[width=0.9\textwidth]{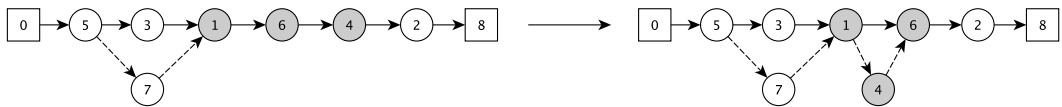}
    \caption{Move Relocate Customer.}
    \label{fig:relocatecustomer}
\end{figure}

This move is based on the method Shift(1,0) described in \cite{Pennaetal2013}. The purpose of this method is to reduce significantly delivery time, therefore removing a customer from the route of the truck and, then, inserting in new drone trip.
First, sub-routes are formed envisioning not to violate the prohibitions.  For example, in Figure \ref{fig:relocatecustomer} the sub-routes are \{0, 5\} and \{1, 6, 4, 2, 8\}. Afterwards, an evaluation is made to determine the launch, the new drone customer, and the returning node. If a triple combination of nodes is feasible, it is then accepted, and the node which is going to be visited by the drone node is removed from the route of the truck and inserted in the route of the drone between the launch and return node. In the example, the customers 1, 4 and 6 were chosen to be the launch, visit and return of the drone, respectively.

\section{Computational Experiments}
\label{section:computational_experiments}

A series of computational experiments was performed to test the effectiveness of the HGVNS. The algorithm was coded in C++ (g++ 5.3.1) and executed on an Intel\textregistered  Core\textsuperscript{TM} i7 Processor 3.6 GHz with 16 GB of RAM running Ubuntu Linux 16.04. The MIP solver used to find a optimal TSP solution was Concorde\footnote{Concorde: http://www.math.uwaterloo.ca/tsp/concorde.html}  3.12.19 with CPLEX 12.6.3. HGVNS was tested in three benchmark sets. Two sets found in literature, advised by \cite{Ponza2016} and by \cite{Agatz2015}. A new set of instances was developed based on the TSPLIB\footnote{TSPLIB: https://www.iwr.uniheidelberg.de/groups/comopt/software/TSPLIB95/} well-known instances for the TSP.

In the tables presented hereafter, the column $Inst. $denotes the instance name and $n$ represents the number of nodes (the depot and $n-1$ customers). Moreover, $BKS$ indicates the best-known solution value found in the literature. 
Columns $s_{FSTSP}$ and $Time$ denote, respectively, the best solution value and the computational time in seconds associated to the MIP running time plus the GVNS algorithm, $gap$ states the difference between $s_{FSTSP}$ and $BKS$. Column $s^*_{TSP}$ describes the optimal tour obtained with the TSP Solver Concorde. Finally, $\overline{{s}_{FSTSP}}$ is the average solution cost of ten runs and $\overline{gap}$ is the gap between the average solutions and the $BKS$. A negative $gap$ indicates an improvement.

\subsection{Results of HGVNS in \cite{Ponza2016} instances}
\label{section:Results_Ponza}

The first set of instances was provided by \cite{Ponza2016} based on the original formulation of \cite{MurrayChu2015}. 
The instances were generated over a map of 32km $\times$ 32km, to reflect the idea of some companies about a 10 km endurance of the drone and thus be able to perform some feasible routes. In the map, the depot is always at coordinate (0, 0) and the locations of the customers are randomly generated. 

The parameters previously defined in this work assumed different values: the service time to the launch ($s^l$) and return ($s^r$) of the drone is 0.6 minutes and 0.5 minutes, respectively. The drone speed is 
80.47 km/h (50 mph), and the truck speed is 56.32 km/h (35 mph). The endurance of the drone is 24 minutes. Additionally, the travel time matrix of both vehicles is calculated based on the same road network, and, finally, the percentage of feasible drone customers is 80\%. 
This experiment was performed by running HGVNS ten times for each instance. Tables~\ref{table:Ponza_tsp}~--~\ref{table:Ponza} provide the results. 

Table~\ref{table:Ponza_tsp} advises the improvement of addressing the FSTSP, i.e., truck and drone in the delivery process, over the classical TSP. According to the results, it is possible to reduce the total travel time up to 30.38\% in this set of instances and get an average improvement of 19.50\%.

\begin{table}[!h]
  \centering
  \caption{FSTSP results versus TSP solutions obtained by Concorde in \cite{Ponza2016} instances.}
    \begin{tabular}{ccccc}
    \hline
     & & & \multicolumn{2}{c}{HGVNS} \\
     \cmidrule(lr){4-5}
    \textit{Inst.} & $n$   & $s^*_{TSP}$ & $gap$ (\%) & $\overline{gap}$ (\%) \\
\hline
    \multicolumn{1}{l}{50.1} & \multicolumn{1}{r}{50} & 14350.60 & -26.79 & -24.34 \\
    \multicolumn{1}{l}{50.2} & \multicolumn{1}{r}{50} & 14458.60 & -24.17 & -23.58 \\
    \multicolumn{1}{l}{50.3} & \multicolumn{1}{r}{50} & 14318.00 & -20.82 & -20.82 \\
    \multicolumn{1}{l}{50.4} & \multicolumn{1}{r}{50} & 14745.10 & -26.37 & -19.10 \\
    \multicolumn{1}{l}{50.5} & \multicolumn{1}{r}{50} & 14564.90 & -28.00 & -24.24 \\
    \multicolumn{1}{l}{100.1} & \multicolumn{1}{r}{100} & 19803.10 & -21.13 & -21.10 \\
    \multicolumn{1}{l}{100.2} & \multicolumn{1}{r}{100} & 19715.30 & -24.43 & -23.27 \\
    \multicolumn{1}{l}{100.3} & \multicolumn{1}{r}{100} & 19644.60 & -26.06 & -26.06 \\
    \multicolumn{1}{l}{100.4} & \multicolumn{1}{r}{100} & 20459.80 & -22.06 & -22.06 \\
    \multicolumn{1}{l}{100.5} & \multicolumn{1}{r}{100} & 19858.20 & -24.72 & -22.05 \\
    \multicolumn{1}{l}{150.1} & \multicolumn{1}{r}{150} & 24519.10 & -27.70 & -26.31 \\
    \multicolumn{1}{l}{150.2} & \multicolumn{1}{r}{150} & 24346.20 & -30.38 & -25.89 \\
    \multicolumn{1}{l}{150.3} & \multicolumn{1}{r}{150} & 25055.50 & -21.64 & -20.00 \\
    \multicolumn{1}{l}{150.4} & \multicolumn{1}{r}{150} & 24038.36 & -15.13 & -2.69 \\
    \multicolumn{1}{l}{150.5} & \multicolumn{1}{r}{150} & 24641.53 & -8.95 & -6.53 \\
    \multicolumn{1}{l}{200.1} & \multicolumn{1}{r}{200} & 27150.20 & -23.84 & -23.84 \\
    \multicolumn{1}{l}{200.2} & \multicolumn{1}{r}{200} & 28303.95 & -2.37 & -2.37 \\
    \multicolumn{1}{l}{200.3} & \multicolumn{1}{r}{200} & 27601.10 & -17.88 & -16.13 \\
    \multicolumn{1}{l}{200.4} & \multicolumn{1}{r}{200} & 28832.90 & -16.35 & -16.35 \\
    \multicolumn{1}{l}{200.5} & \multicolumn{1}{r}{200} & 28458.00 & -25.22 & -23.31 \\ \hline
    \multicolumn{2}{l}{\textbf{Average}}       &       & -21.70 & -19.50 \\
	\hline
    \end{tabular}%
  \label{table:Ponza_tsp}%
\end{table}%

Table~\ref{table:Ponza} compares the FSTSP results provided by \cite{Ponza2016} with HGVNS results.
It is possible to observe that HGVNS acquired better results in all instances, achieving an improvement of 24.84\% in instance \textit{150.2}.
Regarding computational time, HGVNS is executed within a short runtime, 10.15 seconds on average.  However, it is difficult to precisely compare to the method runtime of \cite{Ponza2016} as the work does not report the computer configuration.

\begin{sidewaystable}
  \centering
       \caption{Result running HGVNS in \cite{Ponza2016} instances.}
     \label{table:Ponza}
    \begin{tabular}{ccrrrrrrrrr}
     \hline
     & & & & \multicolumn{2}{c}{SA} & \multicolumn{5}{c}{HGVNS}\\
     & & & &\multicolumn{2}{c}{\cite{Ponza2016}} \\
     \cmidrule(lr){5-6} \cmidrule(lr){7-11}

    Inst. & \multicolumn{1}{c}{\textit{$n$}} &\multicolumn{1}{c}{\textit{$n^\prime$}} & \multicolumn{1}{c}{$BKS$} & \multicolumn{1}{c}{$s_{FSTSP}$} & \multicolumn{1}{c}{\textit{Time (s)}} & \multicolumn{1}{c}{${s}_{FSTSP}$} & \multicolumn{1}{c}{${gap}$ (\%)} & \multicolumn{1}{c}{$\overline{{s}_{FSTSP}}$} & \multicolumn{1}{c}{$\overline{gap}$ (\%)} & \multicolumn{1}{c}{\textit{Time (s)}} \\
  \hline

    \multicolumn{1}{l}{50.1} & \multicolumn{1}{r}{50} & 40    & 12518.93 & 12518.93 & 213.87 & 10506.50 & -16.08 & 10857.02 & -13.28 & 5.57 \\
    \multicolumn{1}{l}{50.2} & \multicolumn{1}{r}{50} & 40    & 12475.14 & 12475.14 & 208.36 & 10964.30 & -12.11 & 11049.21 & -11.43 & 5.91 \\
    \multicolumn{1}{l}{50.3} & \multicolumn{1}{r}{50} & 40    & 12664.65 & 12664.65 & 191.04 & 11336.40 & -10.49 & 11336.40 & -10.49 & 6.39 \\
    \multicolumn{1}{l}{50.4} & \multicolumn{1}{r}{50} & 40    & 12908.18 & 12908.18 & 184.85 & 10856.40 & -15.90 & 11929.50 & -7.58 & 6.67 \\
    \multicolumn{1}{l}{50.5} & \multicolumn{1}{r}{50} & 40    & 12164.83 & 12164.83 & 189.86 & 10486.30 & -13.80 & 11034.30 & -9.29 & 6.76 \\
    \multicolumn{1}{l}{100.1} & \multicolumn{1}{r}{100} & 80    & 17974.85 & 17974.85 & 267.42 & 15618.00 & -13.11 & 15623.84 & -13.08 & 10.70 \\
    \multicolumn{1}{l}{100.2} & \multicolumn{1}{r}{100} & 80    & 17342.18 & 17342.18 & 272.35 & 14899.20 & -14.09 & 15127.50 & -12.77 & 11.75 \\
    \multicolumn{1}{l}{100.3} & \multicolumn{1}{r}{100} & 80    & 17181.88 & 17181.88 & 265.45 & 14524.50 & -15.47 & 14524.50 & -15.47 & 10.06 \\
    \multicolumn{1}{l}{100.4} & \multicolumn{1}{r}{100} & 80    & 18538.03 & 18538.03 & 266.75 & 15947.30 & -13.98 & 15947.30 & -13.98 & 10.62 \\
    \multicolumn{1}{l}{100.5} & \multicolumn{1}{r}{100} & 80    & 17407.43 & 17407.43 & 312.77 & 14948.50 & -14.13 & 15479.22 & -11.08 & 10.08 \\
    \multicolumn{1}{l}{150.1} & \multicolumn{1}{r}{150} & 120   & 22823.38 & 22823.38 & 365.04 & 17728.10 & -22.32 & 18069.32 & -20.83 & 14.05 \\
    \multicolumn{1}{l}{150.2} & \multicolumn{1}{r}{150} & 120   & 22549.55 & 22549.55 & 383.72 & 16949.30 & -24.84 & 18042.20 & -19.99 & 14.11 \\
    \multicolumn{1}{l}{150.3} & \multicolumn{1}{r}{150} & 120   & 23114.14 & 23114.14 & 379.99 & 19633.30 & -15.06 & 20045.00 & -13.28 & 15.62 \\
    \multicolumn{1}{l}{150.4} & \multicolumn{1}{r}{150} & 120   & 22651.00 & 22651.00 & 382.67 & 20400.70 & -9.93 & 23390.90 & 3.27  & 13.33 \\
    \multicolumn{1}{l}{150.5} & \multicolumn{1}{r}{150} & 120   & 22807.41 & 22807.41 & 384.69 & 22435.52 & -1.63 & 23032.05 & 0.98  & 14.58 \\
    \multicolumn{1}{l}{200.1} & \multicolumn{1}{r}{200} & 160   & 26991.21 & 26991.21 & 456.74 & 20676.50 & -23.40 & 20676.50 & -23.40 & 17.38 \\
    \multicolumn{1}{l}{200.2} & \multicolumn{1}{r}{200} & 160   & 27848.14 & 27848.14 & 452.88 & 27632.40 & -0.77 & 27632.40 & -0.77 & 16.41 \\
    \multicolumn{1}{l}{200.3} & \multicolumn{1}{r}{200} & 160   & 27143.78 & 27143.78 & 510.11 & 22665.10 & -16.50 & 23149.00 & -14.72 & 15.89 \\
    \multicolumn{1}{l}{200.4} & \multicolumn{1}{r}{200} & 160   & 28503.18 & 28503.18 & 517.44 & 24119.60 & -15.38 & 24119.60 & -15.38 & 16.25 \\
    \multicolumn{1}{l}{200.5} & \multicolumn{1}{r}{200} & 160   & 27875.87 & 27875.87 & 515.30 & 21281.80 & -23.66 & 21825.40 & -21.71 & 17.20 \\
    \hline
    \multicolumn{2}{l}{\textbf{Average}}       &       &       &       & 336.07 &       & -14.63 &       & -12.21 & 11.97 \\
    \hline
    \end{tabular}
\end{sidewaystable}

\subsection{Results of HGVNS in \cite{Agatz2015} instances}
\label{section:Results_Agatz}

We also address the TSP-D, a problem introduced by \cite{Agatz2015} which presents different restrictions compared to the FSTSP described in this work.
The FSTSP defines an endurance and service time at the launch and return of the drone to the truck.
Thus, some of the constraints of the FSTSP have been relaxed to adapt to the problem. 
First, the drone has unlimited endurance ($e = \infty$), therefore, as result of this definition, the period a vehicle can wait for the other to arrive at the return node is indefinite. Lastly, service time is not required for the launch or return of the drone, i.e., $s^l = s^r = 0$.

The following experiments are using the instances for the TSP-D provided by \cite{Agatz2015}. However, the authors did not provide a complete review of the experiments, thus, we were unable to compare the results obtained with HGVNS and the method used by \cite{Agatz2015}. 

In the original paper, it is explained with details how to generate the coordinates of the locations. The authors proposed three sets of instances changing coordinates distribution.
The first distribution type is called \textit{uniform}. In this instances, the coordinates for every location are drawn independently and uniform. The second type of instance is \textit{single-center} that uses an angle $\alpha$ and a distance $r$ from a normal distribution to obtain instances with locations closer to the center (0, 0). Finally, \textit{double-center} instances use the same strategy of \textit{single-center}, but every location is translated by 200 distance units over the x-axis with a certain probability.

For each type of instance three different scenarios exist, which determines different speed configurations for the drone. A value of  $\alpha = 1$ determines that both vehicles have the same speed. The drone speed is twice as fast as the truck when $\alpha = 2$ and, finally, the drone speed is three times as fast as the truck when $\alpha = 3$.

Tables \ref{table:uniform} -- 
\ref{table:doublecenter} are broken down into the scenarios labeled $\alpha$. The tables present for each group of $n$ the average of ten distinct instances. 
The table's columns follow the same definitions aforementioned. The complete table with the results running HGVNS can be seen in the supplemental material.

\begin{table}[!ht]
\caption{Average solution running HGVNS in the \textit{uniform} instance configuration.}
  \centering
\begin{adjustbox}{max width=\textwidth}
    \begin{tabular}{crrrrrrrrrr}
    \hline
          &       & \multicolumn{3}{c}{$\alpha$ = 1} & \multicolumn{3}{c}{$\alpha$ = 2} & \multicolumn{3}{c}{$\alpha$ = 3} \\
          \cmidrule(lr){3-5}  \cmidrule(lr){6-8}  \cmidrule(lr){9-11}
       
    $n$     & \multicolumn{1}{c}{$s^*_{TSP}$} & \multicolumn{1}{c}{$\overline{{s}_{FSTSP}}$} & \multicolumn{1}{c}{\textit{gap (\%)}} & \multicolumn{1}{c}{\textit{Time (s)}} & \multicolumn{1}{c}{$\overline{{s}_{FSTSP}}$} & \multicolumn{1}{c}{\textit{gap (\%)}} & \multicolumn{1}{c}{\textit{Time (s)}} & \multicolumn{1}{c}{$\overline{{s}_{FSTSP}}$} & \multicolumn{1}{c}{\textit{gap (\%)}} & \multicolumn{1}{c}{\textit{Time (s)}} \\
    \hline
    
    \multicolumn{1}{l}{10} & 7.59  & 6.47  & -14.76 & 0.14  & 5.83  & -23.19 & 0.13  & 5.71  & -24.77 & 0.13 \\
    \multicolumn{1}{l}{20} & 9.80  & 7.93  & -19.13 & 0.11  & 6.84  & -30.20 & 0.85  & 6.86  & -30.05 & 1.15 \\
    \multicolumn{1}{l}{50} & 14.74 & 10.91 & -25.99 & 3.71  & 9.52  & -35.39 & 2.30  & 9.75  & -33.86 & 2.19 \\
    \multicolumn{1}{l}{75} & 17.50 & 12.41 & -29.11 & 16.30 & 10.34 & -40.94 & 10.93 & 10.42 & -40.49 & 10.95 \\
    \multicolumn{1}{l}{100} & 19.52 & 13.98 & -28.39 & 53.50 & 12.03 & -38.36 & 37.77 & 11.96 & -38.72 & 37.26 \\
    \multicolumn{1}{l}{175} & 24.81 & 17.08 & -31.17 & 55.91 & 13.19 & -41.16 & 39.28 & 14.88 & -40.03 & 41.49 \\
    \multicolumn{1}{l}{250} & 29.04 & 19.92 & -31.41 & 185.19 & 17.11 & -41.09 & 191.48 & 17.21 & -40.73 & 189.43 \\ \hline
    \textbf{Average} &       &       & -25.71 & 44.98 &       & -35.76 & 40.39 &       & -35.52 & 40.37 \\
\hline
    
    \end{tabular}%
    \end{adjustbox}
    \label{table:uniform}
\end{table}%

According to the tables it is possible to observe that scenario $\alpha = 1$ obtained the worst results in all set of instances - \textit{uniform}, \textit{single-center} and \textit{double-center} when compared with the other two subcategories of $alpha$. The different solution cost is due to the higher speed of the drone that enables it to visit a greater number of customers, and also visit customers within a further distance. 
The higher improvement occurred in the instance with 75 customers distributed following the configuration \textit{single-center}. The instance achieved an improvement of 62.24\% in comparison to the TSP optimal solution.
Furthermore, the instances presenting a smaller number of customers and the same speed for both vehicles reported the smallest improvement. Additionally, the set $uniform$ showed the lowest average improvement. 
As concerns computational time, the different distribution of customers does not affect runtime, however, when the vehicles present the same speed, and they are uniformly distributed runtime increases. 

\begin{table}[!ht]
  \caption{Average solution running HGVNS in the \textit{single-center} instance configuration.}
    \centering
    \begin{adjustbox}{max width=\textwidth}
    \begin{tabular}{crrrrrrrrrr}
    \hline
          &       & \multicolumn{3}{c}{$\alpha$ = 1} & \multicolumn{3}{c}{$\alpha$ = 2} & \multicolumn{3}{c}{$\alpha$ = 3} \\
          \cmidrule(lr){3-5}  \cmidrule(lr){6-8}  \cmidrule(lr){9-11}

    $n$     & \multicolumn{1}{c}{$s^*_{TSP}$} & \multicolumn{1}{c}{$\overline{{s}_{FSTSP}}$} & \multicolumn{1}{c}{\textit{gap (\%)}} & \multicolumn{1}{c}{\textit{Time (s)}} & \multicolumn{1}{c}{$\overline{{s}_{FSTSP}}$} & \multicolumn{1}{c}{\textit{gap (\%)}} & \multicolumn{1}{c}{\textit{Time (s)}} & \multicolumn{1}{c}{$\overline{{s}_{FSTSP}}$} & \multicolumn{1}{c}{\textit{gap (\%)}} & \multicolumn{1}{c}{\textit{Time (s)}} \\
    \hline

    \multicolumn{1}{l}{10} & 11.27 & 8.87  & -21.25 & 0.14  & 7.76  & -33.77 & 0.14  & 7.81  & -30.73 & 0.13 \\
    \multicolumn{1}{l}{20} & 15.09 & 9.09  & -39.77 & 1.12  & 8.35  & -44.66 & 1.04  & 8.00  & -47.02 & 0.86 \\
    \multicolumn{1}{l}{50} & 21.20 & 13.29 & -37.33 & 3.93  & 8.70  & -58.96 & 2.23  & 8.71  & -58.92 & 2.19 \\
    \multicolumn{1}{l}{75} & 27.80 & 11.87 & -57.31 & 15.32 & 10.50 & -62.24 & 11.18 & 11.38 & -59.04 & 11.39 \\
    \multicolumn{1}{l}{100} & 33.37 & 18.50 & -44.56 & 55.62 & 14.81 & -55.62 & 38.23 & 13.93 & -58.26 & 37.74 \\
    \multicolumn{1}{l}{175} & 45.17 & 23.30 & -48.43 & 53.34 & 17.84 & -60.51 & 43.06 & 17.96 & -60.24 & 41.20 \\
    \multicolumn{1}{l}{250} & 51.00 & 24.30 & -52.35 & 249.92 & 20.45 & -59.91 & 197.12 & 19.63 & -61.51 & 198.72 \\ \hline
    \textbf{Average} &       &       & -43.00 & 54.20 &       & -53.67 & 41.86 &       & -53.67 & 41.75 \\
\hline

    \end{tabular}%
\end{adjustbox}
\label{table:singlecenter}
\end{table}%

\begin{table}[!ht]
  \caption{Average solution running HGVNS in the \textit{double-center} instance configuration.}
    \centering
  \begin{adjustbox}{max width=\textwidth}
    \begin{tabular}{crrrrrrrrrr}
    \hline

          &       & \multicolumn{3}{c}{$\alpha$ = 1} & \multicolumn{3}{c}{$\alpha$ = 2} & \multicolumn{3}{c}{$\alpha$ = 3} \\
          \cmidrule(lr){3-5}  \cmidrule(lr){6-8}  \cmidrule(lr){9-11}

    $n$     & \multicolumn{1}{c}{$s^*_{TSP}$} & \multicolumn{1}{c}{$\overline{{s}_{FSTSP}}$} & \multicolumn{1}{c}{\textit{gap (\%)}} & \multicolumn{1}{c}{\textit{Time (s)}} & \multicolumn{1}{c}{$\overline{{s}_{FSTSP}}$} & \multicolumn{1}{c}{\textit{gap (\%)}} & \multicolumn{1}{c}{\textit{Time (s)}} & \multicolumn{1}{c}{$\overline{{s}_{FSTSP}}$} & \multicolumn{1}{c}{\textit{gap (\%)}} & \multicolumn{1}{c}{\textit{Time (s)}} \\
    
    \hline
    \multicolumn{1}{l}{10} & 18.28 & 13.37 & -26.87 & 0.15  & 12.18 & -33.38 & 0.13  & 11.49 & -37.15 & 0.13 \\
    \multicolumn{1}{l}{20} & 21.25 & 17.27 & -18.71 & 1.39  & 14.48 & -31.87 & 0.98  & 14.34 & -32.52 & 0.99 \\
    \multicolumn{1}{l}{50} & 31.92 & 22.53 & -29.44 & 4.34  & 17.53 & -45.10 & 2.24  & 17.89 & -43.97 & 2.14 \\
    \multicolumn{1}{l}{75} & 39.66 & 25.96 & -34.54 & 19.39 & 20.05 & -49.43 & 11.61 & 19.54 & -50.72 & 11.22 \\
    \multicolumn{1}{l}{100} & 44.18 & 28.91 & -34.56 & 66.39 & 21.24 & -51.92 & 38.22 & 21.82 & -50.62 & 37.86 \\
    \multicolumn{1}{l}{175} & 59.06 & 32.14 & -45.58 & 57.03 & 27.25 & -53.87 & 42.31 & 26.75 & -54.71 & 41.58 \\
    \multicolumn{1}{l}{250} & 69.87 & 38.08 & -45.51 & 260.98 & 31.05 & -55.60 & 193.17 & 30.88 & -55.81 & 196.27 \\ \hline
    \textbf{Average} &       &       & -33.60 & 58.53 &       & -45.88 & 41.24 &       & -46.50 & 41.46 \\ \hline

    \end{tabular}%
    \end{adjustbox}
    \label{table:doublecenter}
\end{table}%

\subsection{Results of HGVNS to new instances from TSPLIB}
\label{section:Results_OurVNS}

Here we propose a new set of instances associated with the original problem specifications described by \cite{MurrayChu2015}. The proposal of new instances set is due to three main reasons. 
First, the instances introduced by \cite{MurrayChu2015} present a small number of customers, up to 10 customers in the FSTSP description. Further, \cite{Ponza2016} considered in his work the same road network for both vehicles, thus to compare our results with the ones of his work, the travel distance was calculated using Euclidean distance. Finally, in \cite{Agatz2015} different restrictions to the FSTSP are stated, such as the endurance of the drone and service time (see Section \ref{section:Results_Agatz}). Therefore, we introduce instances based on the well-known instances of TSPLIB. A total of 25 instances containing between 51 and 200 nodes were selected and appropriately adapted to the problem. 

Drones are usually faster than trucks, since they are not affected by congestion, and they can fly in a straight line not following the street network. Unlike drone, trucks must respect traffic sign regulation and follow the street network. Therefore, it is reasonable to consider different road network for the vehicles. The Euclidean metric is used to describe the drone travel distance considering the straight line flight, and in order to simulate the city block distance, it is used Manhattan distance to represent the truck travel distance.
From the TSP coordinates, two different matrices were calculated to represent the road network that each vehicle travel. The truck matrix matrix was computed using Manhattan distance, and another was estimated using the Euclidean metric to describe the drone travel distance. Additionally, both distances were divided by the speed of each vehicle to obtain the time required to travel among the customers. 

The customers considered eligible to drone delivery are randomly generated such that for every instance there are between 85\% and 90\% serviceable customers and the other 10\% to 15\% are truck-only customers due to geographical limitations, exceedingly heavy parcels or other criteria. 

Furthermore, we considered that both vehicles speed is 40 km/h and the drone endurance is 40 minutes. All these FSTSP characteristics are based on the ones presented by \cite{MurrayChu2015}.

This experiment was performed by running HGVNS ten times for each instance. Table~\ref{table:TSPLIB} provides the results obtained with HGVNS. The table follows the definitions mentioned before. 

HGVNS was able to find better solutions when compared to TSP optimal solutions for all instances.
For instance \textit{pr107}, HGVNS obtained the best solution cost, presenting an improvement of 4\%. Meanwhile, instance \textit{d198} had the minimal improvement with 0.35\%, when compared to the TSP optimal solution value.

\begin{table}[!ht]
  \centering
  \caption{Result running HGVNS in the TSPLIB instances.}
   \centering
   \label{table:TSPLIB}
    \begin{tabular}{crrrrrrrr}
     \hline
      & & & & \multicolumn{5}{c}{HGVNS} \\
      \cmidrule(lr){5-9}

    \multicolumn{1}{c}{\textit{Inst.}} & \multicolumn{1}{c}{$n$} & \multicolumn{1}{c}{$n^\prime$} & \multicolumn{1}{c}{$s^*_{TSP}$} & \multicolumn{1}{c}{$s_{FSTSP}$} & \multicolumn{1}{c}{\textit{gap (\%)}} & \multicolumn{1}{c}{{$\overline{{s}_{FSTSP}}$}} & \multicolumn{1}{c}{$\overline{gap}$ (\%)} & \multicolumn{1}{c}{\textit{Time (s)}} \\
\hline
    \multicolumn{1}{l}{berlin52} & 52    & 46    & 249.75 & 172.25 & -31.03 & 196.25 & -21.42 & 6.50 \\
    \multicolumn{1}{l}{bier127} & 127   & 110   & 3857.00 & 3456.80 & -10.38 & 3587.88 & -6.98 & 23.69 \\
    \multicolumn{1}{l}{ch130} & 130   & 115   & 187.30 & 178.16 & -4.88 & 180.40 & -3.68 & 44.13 \\
    \multicolumn{1}{l}{d198} & 198   & 175   & 463.45 & 461.83 & -0.35 & 461.83 & -0.35 & 27.69 \\
    \multicolumn{1}{l}{eil51} & 51    & 44    & 13.95 & 13.45 & -3.58 & 13.68 & -1.94 & 11.57 \\
    \multicolumn{1}{l}{eil76} & 76    & 65    & 17.50 & 16.35 & -6.57 & 16.68 & -4.69 & 27.14 \\
    \multicolumn{1}{l}{kroA100} & 100   & 86    & 659.90 & 587.80 & -10.93 & 609.71 & -7.61 & 10.95 \\
    \multicolumn{1}{l}{kroA150} & 150   & 134   & 824.25 & 729.95 & -11.44 & 758.43 & -7.99 & 10.95 \\
    \multicolumn{1}{l}{kroA200} & 200   & 180   & 922.45 & 870.65 & -5.62 & 873.99 & -5.25 & 6.53 \\
    \multicolumn{1}{l}{kroB150} & 150   & 133   & 824.80 & 763.15 & -7.47 & 773.72 & -6.19 & 20.20 \\
    \multicolumn{1}{l}{kroB200} & 200   & 178   & 923.00 & 804.47 & -12.84 & 835.89 & -9.44 & 60.00 \\
    \multicolumn{1}{l}{kroC100} & 100   & 87    & 662.70 & 575.30 & -13.19 & 611.31 & -7.75 & 20.07 \\
    \multicolumn{1}{l}{kroD100} & 100   & 90    & 671.20 & 606.45 & -9.65 & 652.34 & -2.81 & 9.15 \\
    \multicolumn{1}{l}{kroE100} & 100   & 86    & 694.65 & 556.30 & -19.92 & 570.30 & -17.90 & 8.57 \\
    \multicolumn{1}{l}{lin105} & 105   & 94    & 662.70 & 378.25 & -42.92 & 380.43 & -42.59 & 10.64 \\
    \multicolumn{1}{l}{pr107} & 107   & 95    & 1866.25 & 1017.50 & -45.48 & 1059.00 & -43.26 & 15.48 \\
    \multicolumn{1}{l}{pr124} & 124   & 107   & 2048.50 & 1553.80 & -24.15 & 1566.62 & -23.52 & 25.45 \\
    \multicolumn{1}{l}{pr136} & 136   & 119   & 3015.25 & 2542.00 & -15.70 & 2559.00 & -15.13 & 44.50 \\
    \multicolumn{1}{l}{pr144} & 144   & 126   & 2055.00 & 1666.25 & -18.92 & 1675.75 & -18.45 & 62.63 \\
    \multicolumn{1}{l}{pr152} & 152   & 132   & 2321.45 & 1919.35 & -17.32 & 1938.47 & -16.50 & 70.33 \\
    \multicolumn{1}{l}{rat99} & 99    & 87    & 37.95 & 37.15 & -2.11 & 37.33 & -1.63 & 9.76 \\
    \multicolumn{1}{l}{rat195} & 195   & 171   & 72.25 & 71.40 & -1.18 & 71.64 & -0.84 & 24.13 \\
    \multicolumn{1}{l}{rd100} & 100   & 86    & 248.93 & 240.46 & -3.40 & 243.84 & -2.04 & 10.83 \\
    \multicolumn{1}{l}{st70} & 70    & 63    & 21.50 & 20.50 & -4.65 & 21.00 & -2.33 & 3.85 \\ \hline
    \textbf{Average} &       &       &       &       & -13.49 &       & -11.26 & 23.53 \\ \hline
    \end{tabular}%
  \label{tab:addlabel}%
\end{table}%

\begin{figure}[!ht]
\centering
   \includegraphics[width=0.5\linewidth, angle=-90]{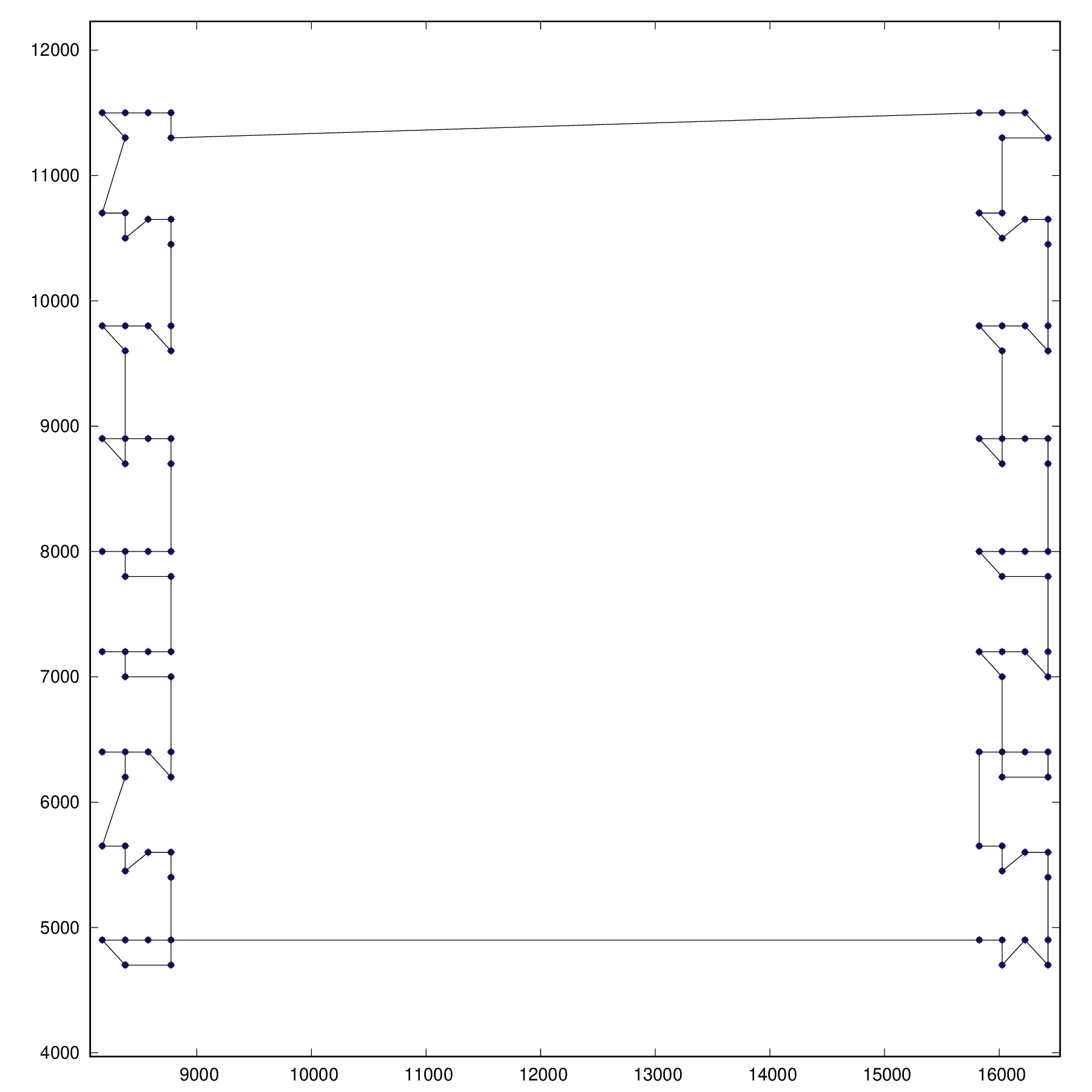}
   \caption{Delivery route with truck only.}
   \label{fig:noUAV} 
\end{figure}

\begin{figure}[!ht]
\centering
    \includegraphics[width=0.5\linewidth, angle=-90]{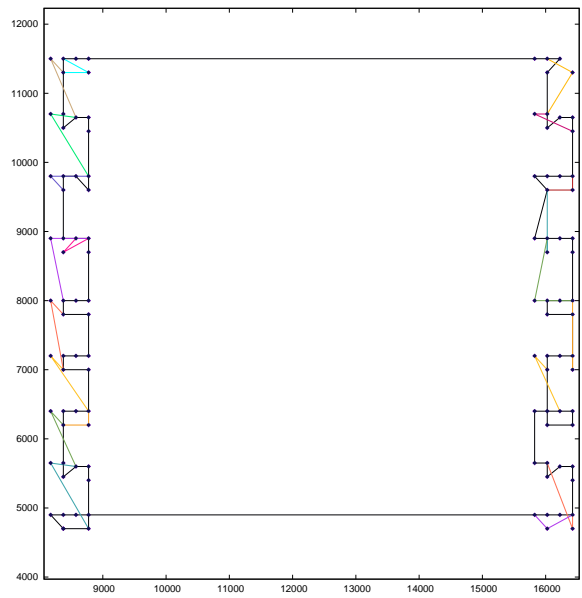}
   \caption{Delivery route with truck and drone.}
   \label{fig:withUAV}
\end{figure}

Figures \ref{fig:noUAV} and \ref{fig:withUAV} shows the solution for \textit{pr107} instance, where the colored lines represent the UAV routes while the continuous black line represents the truck route. HGVNS manages to assign 16 delivers to the drone. It is possible to observe that the triangulations performed by the truck in the route illustrated by Figure \ref{fig:noUAV} were, in its majority, replaced by a drone trip illustrated by Figure \ref{fig:withUAV}, transforming the truck trip into a more orthogonal one.

\section{Concluding Remarks}
\label{section:conclusions}

The Flying Sidekick Traveling Salesman Problem (FSTSP) concerns a variant of TSP which has been showing potential over the last years through the constant investments of companies such as JD.com, Amazon, Mercedes-Benz, among others.
The problem consists in the use of Unmanned Aerial Vehicles (UAV), also known as drones, working collaboratively with trucks in parcel delivery.

We presented a hybrid heuristic using the complementary characteristics of truck and drone to perform deliveries with reduced time. The algorithm implemented in this work, named HGVNS, initially uses a MIP solver to obtain the optimal TSP tour, whose solution is subsequently enhanced by the meta-heuristic General Variable Neighborhood Search (GVNS).

HGVNS was tested in three benchmark sets, two sets found in literature: the instances presented by \cite{Ponza2016} and the instances introduced by \cite{Agatz2015}. The third set is a new one developed based on TSPLIB instances for the TSP. 

In the instances introduced by \cite{Ponza2016}, HGVNS improved the solution of the majority instances achieving an improvement up to 24.84\% to the previously best-known solution values.

A variant of the FSTSP proposed by \cite{Agatz2015} called TSP-D was also studied. The computational tests using their 1383 instances evidenced that the speed of the drone interferes with the total delivery time, however, the effect of double or triple truck speed is very similar. That said, the best improvement occurred in one instance with 75 customers where the drone traveled twice as fast as the truck, and the least improvements were observed for instances where both vehicles presented the same speed.

A new set of instances based on the well-known instances of TSPLIB is proposed to fulfill the need large instances following the original model of \cite{MurrayChu2015}. The best solution found in these instances shows an improvement of 45.48\% when over the optimal TSP tour value.

A new modality of parcel distribution is raising from the increasing development of drones and the effort of companies to perform deliveries faster at a reduced cost. Thus, this work has plenty to contribute demonstrating that collaborative work of truck and drone can drastically decrease delivery times up to~67.79\%. 

Future research directions include formulating a Mixed Integer Linear Program for the Flying Sidekick Traveling Salesman Problem. Moreover, it opens a huge field of research in distribution and logistics area. For example, one can study the capacitated version of the problem with multiple delivery trucks and multiple drones per truck.

\section*{Acknowledgments}
This research was supported by the Foundation for Research Support of the State of Minas Gerais (FAPEMIG).


\nocite{*}

\bibliographystyle{itor}
\bibliography{itor}

\end{document}